\documentclass[10pt,reqno]{amsart}
\usepackage{amsmath,bm}
\usepackage{amsfonts, fourier}
\usepackage{amssymb}
\usepackage{amsthm}
\usepackage{amscd}
\usepackage{algorithm}
\usepackage[noend]{algpseudocode}
\makeatletter
\def\BState{\State\hskip-\ALG@thistlm}
\makeatother
\usepackage{amsmath,tikz}
\usepackage{amssymb}
\usepackage{amsthm}
\usepackage{graphics}
\usepackage{eucal}
\usepackage{mathrsfs}

\usepackage{blkarray}

\usepackage[utf8]{inputenc} 
\usepackage{textcomp} 
\usepackage{graphicx}  
\usepackage{flafter}  
\usepackage{bm}  
\usepackage{tikz}
\usepackage[pdftex,bookmarks,colorlinks,breaklinks]{hyperref}  
\usepackage{color}
\usepackage{memhfixc}  
\usepackage{pdfsync}  
\usepackage{amssymb}
\usepackage{amsfonts}
\theoremstyle{plain}
\newtheorem*{maintheorem*}{Main Theorem}
\newtheorem*{thm*}{Theorem}
\newtheorem*{thma*}{Theorem A}
\newtheorem*{thmaa*}{Theorem A'}
\newtheorem*{thmb*}{Theorem B}
\newtheorem*{thmo*}{Theorem 1.1}
\newtheorem*{thmc*}{Theorem C}
\newtheorem*{thmd*}{Theorem D}
\newtheorem*{thmf*}{Theorem 4.1}
\newtheorem*{remark*}{Remark}
\newtheorem*{conjecture*}{Conjecture}
\newtheorem*{prop*}{Proposition}
\newtheorem*{lem*}{Basic Lemma}
\newtheorem{thm}{Theorem}[section]

\newtheorem{lem}[thm]{Lemma}
\newtheorem{prop}[thm]{Proposition}

\newtheorem*{proofc*}{Proof of Theorem C}

 \usetikzlibrary{calc}

\def\sinc{~\mathrm{sinc}}

\begin{document}

\author{Youssef Lazar}
\email{ylazar77@gmail.com}

\title{
A lattice point counting approach for the study of the number of self-avoiding walks on  $\mathbb{Z}^{d}$}
  
\maketitle

\begin{abstract}    We reduce the problem of counting self-avoiding walks in the square lattice to a problem of counting the number of integral points in multidimensional domains. We obtain an asymptotic estimate of the number of self-avoiding walks of length $n$ in the square lattice.  This new formalism gives a natural and unified setting in order to study the properties of the number of self-avoiding walks in the lattice $\mathbb{Z}^{d}$ of any dimension $d\geq 2$. 
\end{abstract}


\section{Introduction}

\noindent \subsection{The origins of the problem} Self-avoiding walks has been introduced by Orr \cite{orr} and Flory \cite{flory} around 1940 in order to model a long chain polymer in a three-dimensional lattice.   If we denote by $\Omega$ the space of all possible configurations, then,  for each configuration $\omega \in \Omega$ one can assign an weight $H(w)$ to the corresponding configuration. This weigth function $H$ can be used to define an energy functional $H$, namely the \textit{Hamiltonian }of the system. For instance, in the \textit{Ising model} an energy function is given by 
 the number of neighbours which differ from each other. 
 Boltzmann's thermodynamical formalism for the ferromagnetism gives the probability of a configuration to occur at a fixed temperature $T$  when the space of configuration is given by spins $\Omega=\{ \pm 1\}^{n}$, 
$$ \mathbf{P}[\omega] = \dfrac{e^{-\beta H(\omega)}}{Z_{\beta}}$$
where $\beta=1/T$ given  the normalization $Z_{\beta} = \sum_{\omega \in \Omega} e^{-\beta H(\omega)}.$
This probability has a sense only when the partition function satisfies  $Z_{\beta} < \infty$. It is not diffcult to see that there exists a critical exponent $\beta_{c}$ such that $Z_{\beta}$ is finite if $\beta > \beta_{c}$ and infinite otherwise. Also, one notes that $H(\omega)=n$ if and only if $\omega$ defines a self-avoiding path of length $n$ for the Ising model. Hence, 
$$ Z_{\beta} = \sum_{n \geq 1} c_{n} e^{-\beta n}.$$

\noindent  The partition function for the Ising model involves the coefficients $c_{n}$ which corresponds to the number of self-avoiding walks of length exactly $n$ in the square lattice $\mathbb{Z}^{2}$.  In particular, the knowledge of $c_{n}$ would give the critical exponent  $\beta_c $ for the Ising model. More precisely, if we denote $\mu_c = e^{\beta_{c}}$, then it corresponds to the radius of convergence of the power series
$$Z(x)= \sum_{n \geq 0} c_{n} x^{n}.$$

\noindent An argument of sub-multiplicativty due to Hammersley  shows that $c_{n}$ has an exponent growth in the sense that 
$$ \lim_{n} c_{n}^{1/n} = \mu_c.$$
The number $\mu_c$ is a universal constant called the \textit{connective constant} of the square lattice. The exact value of $\mu_{c}$ is still unknown for the square lattice but it is not difficult to show that $2 <\mu_{c} < 3$ and that necessarily $\mu^{n}_{c} \leq c_{n}$ for all $n \geq 1$. Hammersley and Welsh \cite{hw} gave the following best known upper bounds 
$$ c_{n} \leq \mu_{c}^{n} e^{c \sqrt{n}}.$$
In the early eighties, the physicist B. Nienhuis \cite{nienhuis} has conjectured that the expected exact asymptotic behaviour of $c_{n}$ as the length tends to infinity should be
$$ c_{n}  \sim A n^{11/32} \mu_{c}^{n}$$
for some universal positive constant $A$. One can note that, while the connective constant $\mu$ depends essentially on the lattice, the hypothetical exponent $11/32$ is universal to any lattice of the plane.  There is a highly nontrivial case when the connective constant is known explicitly which is the case of the hexagonal lattice in the plane. Indeed, quite recently, Duminil-Copin and Smirnov \cite{dmhex} found the exact value of the connective constant for the hexagonal lattice predicted by Nienhuis again, namely $$\mu_{hex.} = \sqrt{2+\sqrt{2}}.$$

\noindent A result due to Kesten \cite{kesten} tells us that 
$$ \lim_{n} \dfrac{c_{n+2}}{c_{n}} = \mu_{c}^{2}$$
for any lattice. Kesten's result is coherent with Nienhuis prediction. However, a proof of the existence of the limit $\frac{c_{n+1}}{c_{n}} =\mu_{c}$ is still an open problem. The difficulty of the enumeration of self avoiding walks is due to the fact all moves are correlated to each other.  The self-avoidance feature gives rise to trapped walks which are of length $n$ and which cannot be extended to a self-avoiding walk of length $n+1$.  \\

From a probabilistic point of view, the scaling limit of self-avoiding walk as the length tends to infinity is predicted to converging towards a continuous conformally invariant random process called the \textit{Schramm-Loewner evolution } \cite{sle}.  This process has been studied intensively during the last twenty years and lot is known about them. In particular, using this conjectural limit distribution i.e. the $SLE(8/3)$, one can confirm the prediction of Nienhuis. For classical random walks , it  correspond to Donsker's theorem which asserts that the scaling limit of a random walk is the Brownian motion. \\

 For a typical self-avoiding walk of length $n$, it is natural to expect that its endpoint does go further from the origin than a typical random walk. For a random walk of length $n$, with self-intersections allowed, its mean displacement of the endpoints from the origin is of order $\sqrt{n}$, whereas for self-avoiding walks,  it is conjectured to be $n^{3/4}$ which is much higher, this bound has been predicted by the chemist Paul Flory \cite{flory}.  For a complete survey of self-avoiding walks, we refer the reader to \cite{bauer}.

\subsection{Set-up and goal of the paper} A walk of length $n$ in the lattice $\mathbb{Z}^{2}$ can be encoded by a set of $n$ successive points $w_{0},$$ \ldots $, $ w_{n-1}, w_{n}$ in $\mathbb{Z}^{2}$ such that $\|w_{i+1}-w_{i} \|=1$ for $i=0, \ldots, n-1$. Each point $w_{i}$ correspond to the position of the walk at step $i$. All the walks we are going to consider starts at the origin i.e. $w_{0}=0$. Thus, a walk of length $n$, say $w=(w_{1}, \ldots, w_{n})$ can be seen as an element of $\mathbb{Z}^{2n}$.  \\
A self-avoiding walk of length $n$ in the square lattice is a walk $\omega=(w_{0}, \ldots, w_{n})$ which has no self-intersections i.e. $w_{j} \neq w_{k}$ for all $1 \leqslant j < k \leqslant n$. 
We denote by $\Omega_{n}$,  the set of all self-avoiding walks of length $n$ in the square lattice and $c_{n}$ is the number of self-avoiding walks of length $n$ in the square lattice, i .e. $c_{n} = |\Omega_{n}|$. The number of lattice points in any Borel measurable bounded domain $B$ of $\mathbb{R}^{n}$ is intimately linked to the discrete uniform distribution on a finite Borel subset $E$ of $\mathbb{R}^{n}$, which is defined by 
$$ \mathbf{P}(X \in B) = \dfrac{|B \cap \mathbb{Z}^{n}|}{|E|} = \dfrac{1}{|E|} \sum_{x\in \mathbb{Z}^{n} } \mathbf{1}_{B}(x).$$

\noindent In order to obtain a closed formula for the discrete uniform probability distribution $ \mathbf{P}_{n}$ on the set $\Omega_{n}$ of self-avoiding walks of length $n$, one has to compute $c_{n}=|\Omega_{n}|$, hence for any Borel subset $A$ in $\Omega_{n}$, 
$$ \mathbf{P}_{n}(\omega \in A) = \dfrac{| A |}{|\Omega_{n}|} = \dfrac{1}{c_{n}} \sum_{ \omega \in \Omega_{n} } \mathbf{1}_{A}(\omega).$$

\noindent The aim of the paper is to study $c_{n}$ from a purely theoretical point of view, in that, we are going to translate the problem of estimating $c_n$ to a problem of counting lattice points in some complicated domain in dimension increasing with $n$.  Several methods are available to do that, but the most efficient is the use of Fourier transforms. Indeed, this leads quite rapidly to some non-trivial estimate for the number of self-avoiding walks in any dimension by using a unified formalism. 
\noindent  The main idea of the paper relies on a key formula for $c_{n}$ which we will obtain  in (\ref{uniformstyle}), 
$$ c_{n}= \sum_{\sigma \in \mathbb{Z}^{2n}}\prod_{0 \leqslant j < k \leqslant n} \left( \chi_{B^{+}_{jk}(\sigma)} \left( \sigma_{k}, \sigma_{k+n} \right) \prod_{l\neq j,k} \chi_{[-l,l]^{2}}(\sigma_{l}, \sigma_{l+n}) \right).$$
This formula gives an arithmetical meaning to $c_{n}$, in that, it counts the number of integral vectors in $\mathbb{Z}^{2n}$ which lie in some complicated bounded domain of $\mathbb{R}^{2n}$. Roughly speaking, to each self-avoiding walk of length $n$ we can associate a vector $\sigma=(\sigma_{1}, \ldots, \sigma_{n}; \sigma_{n+1}, \ldots, \sigma_{2n}) \in \mathbb{Z}^{2n}$ such that 
for any $0 \leq j < k \leq n $, the two-dimensional vector $(\sigma_{k}, \sigma_{k+n}) $  lies in the intersection $[-k,k]^{2} \cap \left(  B((\sigma_{j}, \sigma_{j+n}), (k-j)^{2})\backslash B((\sigma_{j}, \sigma_{j+n}), 1/2) \right)$ where we set $\sigma_{0}=0$.
Combining this fact with Poisson's summation formula, we are able to derive an asymptotic estimate not only for $c_{n}$, but also for the number of self-avoiding walks of length in $\mathbb{Z}^{n}$ for any dimension $d \geq 1$.

\subsection{Main results}

An asymptotic estimate is obtained for number of self avoiding walks in the square lattice by using a lattice point counting in certain multidimensional domains.

\begin{thm}\label{main} As $n$ gets large,  the number  $c_n$ of self-avoiding walks of length $n$ behaves as follows, 
\begin{equation}\label{formulation}
c_{n} \sim \sum_{v \in \mathbb{Z}^{2n}} \int_{\mathbb{R}^{2n}} \ldots \int_{\mathbb{R}^{2n}} \widehat{\psi_{1,2}}(\xi_{1,2}) \ldots \widehat{\psi_{n-2,n}}(\xi_{n-2,n}) \widehat{\psi_{n-1,n}}\left( v -\sum_{1\leqslant j<k\leqslant n-1 }\xi_{j,k} - \sum_{1 \leqslant j \leqslant n-2} \xi_{j,n}\right) \left(\bigotimes_{1 \leqslant j<k\leqslant n-1} d\xi_{j,k}\right) \bigotimes \left(\bigotimes_{1\leqslant j \leqslant n-2} d\xi_{j,n} \right).
 \end{equation}
 where for each $\xi \in \mathbb{R}^{2n}$, 
 \begin{equation}\label{fourierfunctions}
\begin{array}{c} \widehat{\psi_{jk}}(\xi) = \left(\displaystyle \prod_{{1\leqslant l \leqslant n}\atop{l \neq j,k}} 
l^{2}\sinc (2\pi l \xi_{l}) \sinc(2\pi l \xi_{l+n}) \right)  j^{2}\sinc(2\pi j (\xi_{j}+\xi_{k})) \sinc(2\pi j (\xi_{j+n}+\xi_{k+n})) \\
\left\lbrace \dfrac{(k-j) J_{1}(2\pi (k-j)\|(\xi_{k}, \xi_{k+n}) \|)  - \frac{1}{2} J_{1}(\pi \|(\xi_{k}, \xi_{k+n}) \|)  }{\|(\xi_{k}, \xi_{k+n}) \|} \right\rbrace\end{array}
\end{equation}
and $J_{1}$ is the first order Bessel function of the first kind. 
 \end{thm}
\noindent \textbf{Remark.} The main term of the asymptotic is achived for by the term $v=0$ in the sum (the volume term in Poisson Summation Formula) due to the fact that the Fourier coefficients decay to zero (Riemann-Lebesgue).

Applying exactly the same strategy, we are able to derive without extra effort, the following estimate on the number of SAWs in any dimension $d \geq 2$.

\begin{thm}\label{maind} As $n$ gets large, the number $c_n(d)$  of SAWs of length $n$ in $\mathbb{Z}^{d}$ behaves as follows, 
$$c_{n}(d) \sim \sum_{v \in \mathbb{Z}^{dn}} \int_{\mathbb{R}^{dn}} \ldots \int_{\mathbb{R}^{dn}} \widehat{\psi_{1,2}}(\xi_{1,2}) \ldots \widehat{\psi_{n-2,n}}(\xi_{n-2,n})$$ $$ \widehat{\psi_{n-1,n}}\left( v -(\sum_{1\leqslant j<k\leqslant n-1 }\xi_{j,k} + \sum_{1 \leqslant j \leqslant n-2} \xi_{j,n})\right) \left(\bigotimes_{1 \leqslant j<k\leqslant n-1} d\xi_{j,k}\right) \bigotimes \left(\bigotimes_{1\leqslant j \leqslant n-2} d\xi_{j,n} \right)$$
where the Fourier transform  for each $\xi \in \mathbb{R}^{dn}$ is given by
  \begin{equation}\label{fourierfunctions}
\begin{array}{c} \widehat{\psi_{jk}}(\xi) = \left(\displaystyle \prod_{{1\leqslant l \leqslant n}\atop{l \neq j,k}} 
l^{2}\sinc (2\pi l \xi_{l}) \sinc(2\pi l \xi_{l+n})  \ldots  \sinc(2\pi l \xi_{l+dn}) \right)  j^{2}\sinc(2\pi j (\xi_{j}+\xi_{k})) \ldots  \sinc(2\pi j (\xi_{j+dn}+\xi_{k+dn}))\\
\left\{ \dfrac{(k-j) J_{d/2}(2\pi (k-j)\|((\xi_{k}, \xi_{k+n},  \ldots, \xi_{k+dn}) \|)  - \frac{1}{2} J_{d/2}(\pi \|(\xi_{k}, \xi_{k+n},  \ldots, \xi_{k+dn}) \|)  }{\|(\xi_{k}, \xi_{k+n}, \ldots, \xi_{k+dn} ) \|^{d/2}} \right\}
\end{array}
\end{equation}
and $J_{d/2}$ is the Bessel function of the first kind of order $d/2$. 
 \end{thm}
 
\newpage
\section{Reduction to a lattice point counting in multidimensional domains}~

\subsection{Arithmetical interpretation of $c_{n}$} \noindent  We start by the natural reformulation of $c_{n}$ using indicator functions, 
\begin{equation}\label{defcn}
c_n = \sum_{ w \in \mathbb{Z}^{2n}} \prod_{0 \leqslant j < k \leqslant n} \chi \left(1 \leqslant \| w_{k}-w_{j} \| \leqslant k-j \right).
\end{equation}

\noindent Here $w=(w_{1}, w_{2}, \ldots,w_{n}) \in \mathbb{Z}^{2n}$ where for each $1 \leqslant i \leqslant n$, $w_{i}$ is the location of the walk after the \textit{ith}-move in the lattice $\mathbb{Z}^{2}$ and $w_{0}=(0,0)$ is the initial point of the walk. The formula (\ref{defcn}) is coherent since the right hand side counts the number of walks of length which satisfies that $$\| w_{k}-w_{j} \| \geq 1$$ which is equivalent to asking $w_{k}\neq w_{j}$ since the norm of integral vectors assumes its value in the set of nonnegative integers. Concerning the reverse inequalities 
$$  \| w_{k}-w_{j} \| \leqslant k-j ~, ~~~0 \leqslant j<k \leqslant n $$ they come from the fact that if we denote the \textit{ith}-move by by $v_{i}=(x_{i}, y_{i}) \in \{(\pm1,0), (0,\pm 1) \}$ we get that 
$$ w_{i}= v_{1}+ \ldots + v_{n}= (x_{1}+ \ldots+ x_{i}, y_{1}+ \ldots+y_{i})$$ with the orthogonality condition $x_{i}y_{i}=0$.
In particular, this implies that,
\begin{equation}
\|w_{k}-w_{j}\|^{2} = (x_{j+1}+ \ldots+ x_{k})^{2} +  (y_{j+1}+ \ldots+ y_{k})^{2}
\end{equation}
 is maximized by the value $(k-j)^{2}$\footnote{This maximum is obtained when the walk realizes the same move between the $jth$ step and the $kth$ step for instance $v_{j+1}= \ldots = v_{k}=e_{1}$ in this case,  $x_{j+1}+ \ldots+ x_{k} =k-j$ and $ y_{j+1}+ \ldots+ y_{k}=0$ } by orthogonality (exactly one only between $x_{i}$ and $y_{i}$ has to be nonzero since $x_{i}^{2}+ y_{i}^{2}=1$ i.e. 
$\| w_{i+1}- w_{i} \|=1$ for each $i=0, \ldots, n-1$). This shows that formula (\ref{defcn}) counts the number of self-avoiding walks of length in the square lattice. \\ 
For pratical reasons, we denote by $v_{i}=(x_{i}, x_{n+i})$ the \textit{ith}-move of a self avoiding walk of length $n$. Hence, a self avoiding walk $w \in \Omega_{n}$ is entirely determined by the \textit{total displacement vector} $$x=(x_{1},x_{2}, \ldots, x_{n}; x_{n+1}, \ldots, x_{2n}) \in \{0,\pm 1 \}^{2n}.$$ Therefore, the formula (\ref{defcn}) can be rewritten as follows,

  \begin{equation}\label{defcnx}
c_n = \sum_{{x \in  \{0,\pm 1 \}^{2n}}} \prod_{0 \leqslant j < k \leqslant n}\chi \left(1 \leqslant  (x_{j+1}+ \ldots+ x_{k})^{2} +  (x_{j+n+1}+ \ldots+ x_{k+n})^{2} \leqslant (k-j)^{2} \right).
\end{equation}

\noindent In order, to apply Poisson summation formula, we use the obvious fact that $x \in  \{0,\pm 1 \}^{2n}$ if and only if $x \in \mathbb{Z}^{2n}$ and $\|x\|_{\infty} \leqslant 1$, therefore 
 \begin{equation}\label{defcnx}
c_n = \sum_{{x \in \mathbb{Z}^{2n}}\atop{\|x\|_{\infty} \leqslant 1}} \prod_{0 \leqslant j < k \leqslant n}\chi \left(1 \leqslant  (x_{j+1}+ \ldots+ x_{k})^{2} +  (x_{j+n+1}+ \ldots+ x_{k+n})^{2} \leqslant (k-j)^{2} \right).
\end{equation}

\subsection{The geometry of the problem}~~\\

\noindent  \textit{Geometrical observation: } The previous formulation (\ref{defcnx}) shows that finding $c_{n}$ reduces the problem to counting the number of lattice points in the unit hypercube $[-1,1]^{2n}$ cuted out by ellipsoids shells of dimension $2(k-j)$ defined by the inequalities

\begin{equation}
1 \leqslant  (x_{j+1}+ \ldots+ x_{k})^{2} +  (x_{j+n+1}+ \ldots+ x_{k+n})^{2} \leqslant (k-j)^{2}.
\end{equation}

 \noindent For $0 \leqslant j < k \leq n$, let us introduce the set 
$$\mathcal{E}_{jk} = \{ v \in [-1,1]^{2n},  1 \leqslant  (x_{j+1}+ \ldots+ x_{k})^{2} +  (y_{j+1}+ \ldots+ y_{k})^{2} \leqslant (k-j)^{2}  \}. $$
It is the intersection of the hypercube $ [-1,1]^{2n}$ with the ellipsoid shell of inner radius 1 and outer radius $(k-j)^{2}$.  These ellispoids are defined by the quadratic forms $$Q_{jk}(x) =  (x_{j+1}+ \ldots+ x_{k})^{2} +  (y_{j+1}+ \ldots+ y_{k})^{2}$$ which are clearly degenerated and isotropic. The isotropic cone of $Q_{jk}$ is the intersection of two hyperplanes with respective normal vectors $n_{jk}= e_{j+1}+ \ldots+ e_{k}$ and $n^{\prime}_{jk}=e_{j+n+1}+ \ldots+ e_{k+n}$ in $\mathbb{R}^{2n}$. Hence, 
$$\mathcal{E}_{jk} = \{ 1 \leq Q_{jk} \leq k-j\}$$

\noindent We continue by defining $\mathcal{E}_{n}$ as the intersection of the $\mathcal{E}_{jk}$ ($0 \leqslant j < k \leqslant n $), i.e.
$$ \mathcal{E}_{n} = \bigcap_{0 \leqslant j < k \leqslant n} \mathcal{E}_{jk}.$$

\noindent We remark immediately that $\Omega_{n}= \mathcal{E}_{n} \cap \mathbb{Z}^{2n}$ and the crucial fact that $c_{n}$ is the number of lattice points of the set $\mathcal{E}_{n}$. The task of counting lattice points in a given domain is a difficult problem in general.  A class of domain for which the problem can be solved is the case of centrally symmetric convex bodies. In that case, the number of lattice points is proportional to the volume of $\mathcal{E}_{n}$, so, we are reduced to a volume computation which is in general more approachable.  \\
In our situation,  we can show that $\mathcal{E}_{n}$ is a convex body and we are reduced to  the computation of the volume of $\mathcal{E}_{n}$.  We are only interested in the asymptotic of  the volume of $\mathcal{E}_{n}$ as $n \rightarrow \infty$, so we can avoid a closed formula for the volume and thus, for $c_{n}$. \\

\subsection{Poisson Summation Formula}  Poisson summation formula asserts that for a function in the Schwartz class (of rapidly decreasing functions) one has the equality
$$ \sum_{x \in  \mathbb{Z}^{n}} f(n) = \sum_{\xi \in  \mathbb{Z}^{n}} \widehat{f}(\xi).$$
 
\noindent  The classical strategy to find the number of lattice points of $\mathcal{E}_{n}$ is to proceed as follows, if $f=\chi_{\mathcal{E}_{n}}$  the indicator function for $\mathcal{E}_{n}$ then using Poisson formula, one formally has 
 
 $$ \sum_{x \in  \mathbb{Z}^{n}} \chi_{\mathcal{E}_{n}}(n) = \sum_{\xi \in  \mathbb{Z}^{n}} \widehat{\chi_{\mathcal{E}_{n}}}(\xi).$$

\noindent The left hand side coincides with the number of lattice points in our compact $\mathcal{E}_{n}$, hence we obtain a formula for $c_{n}$. The only thing we need to precise is that $f=\chi_{\mathcal{E}_{n}}$ has to be regularized in be in $\mathcal{S}$.  This does not affect the asymptotics. Finally, we get the following approximation for the number of self-avoiding walk of length $n$
$$ c_{n}=\sum_{x \in  \mathbb{Z}^{n}} \chi_{\mathcal{E}_{n}}(n)  \sim   \sum_{\xi \in  \mathbb{Z}^{n}} \widehat{\chi_{\mathcal{E}_{n}}}(\xi).$$
Hence our strategy will be to estimate the RHS as $n \rightarrow \infty$. In fact, we do not need to estimate the whole sum since we are only interested by the main term of the asymptotic not with the error term. We know that the Fourier coefficients decreases to zero, hence the main term is the volume term 
$$\widehat{\chi_{\mathcal{E}_{n}}} (0) = \int_{\mathbb{R}^{2n}} \chi_{\mathcal{E}_{n}}(x)  dx = \mathrm{vol} (\mathcal{E}_{n}). $$
As  we have noted before we are not able to compute the volume of $\mathcal{E}_{n}$ directly. Instead, we will compute the coefficient $\widehat{\chi_{\mathcal{E}_{n}}} (0)$
Finally, we get 
$$ c_{n} \sim \widehat{\chi_{\mathcal{E}_{n}}} (0) (1+o(1)).$$

\noindent In the present study, we are going to find a nice integral formulation of $\widehat{\chi_{\mathcal{E}_{n}}} (0)$  involving products of explicit transcendental functions. The key Lemma \ref{cnsigma} give a more suitable form for the counting by using of a unimodular change of variables $x \leftrightarrow \sigma$  which makes things easier (without this change of variables, it almost impossible to handle the factors arising from Fourier transformation). We discuss these points in full detail in the following sections.

\section{Proof of Theorem \ref{main}}
\noindent We begin by the following Lemma which gives more appreciable formulation of $c_{n}$.
\begin{lem}\label{cnsigma}
$$c_n = \sum_{{\sigma \in \mathbb{Z}^{2n}}\atop{|\sigma_{i~\mathrm{mod}(n)}| \leqslant i}} \prod_{0 \leqslant i < k \leqslant n}\chi \left(1/2  \leqslant  (\sigma_{k}-\sigma_{i})^{2} + (\sigma_{k+n}-\sigma_{i+n})^{2} \leqslant( k-i)^{2} \right).$$
\end{lem}
\noindent \textbf{Proof.}  Let $x=( x_{1}, \ldots x_{n} , x_{n+1},  \ldots, x_{2n})\in \mathbb{Z}^{2n}$ with $\| x\|_{\infty} \leq 1$, this simply means that $x_{i}, x_{i+n} \in \{0, \pm 1 \}$ for each $1 \leqslant i \leqslant 2n$. For every $1 \leqslant i \leqslant n$, let us introduce the couple of partial sums of order $(i,i+n)$ 
\begin{center}
$\sigma_{i}= x_{1}+ \ldots + x_{i}$ \ \ \ for $1 \leqslant i \leqslant 2n$.
\end{center}
From the fact that $x_{i}x_{i+n} =0$  for all $1 \leqslant i \leqslant n$,  we obtain that, 
$$ (\sigma_{i}, \sigma_{i+n}) \in [-i,i]^{2}.$$
In particular, we have the matrix relation
 $$  \left[  \begin{array}{c}
\sigma_{1}  \\
\sigma_{2} \\
 \vdots \\
 \sigma_{n-1}  \\
 \sigma_{n} \\\hline
 \sigma_{n+1}\\
 \sigma_{n+2}\\
 \vdots \\
 \sigma_{2n-1}\\
 \sigma_{2n}
\end{array} \right]= \left[ \begin{array}{@{}cccccc|ccccc@{}}
    1 & 0 & \ldots & 0 & 0 &0 &0  &\ldots &0 &0 &0 \\
    1 & 1 & \ldots & 0 & 0 &0 & 0&   \ldots &0  &0  & 0 \\
    \vdots & \vdots &\vdots & \vdots & \vdots & \vdots & \vdots & \vdots & \vdots & \vdots & \vdots\\
     1 & 1 &  1 & \ldots  & 1& 0& 0 & 0 & 0 & \ldots  & 0 \\
         1 & 1 & \ldots & 1 & 1 &1 & 0& 0 &  \ldots &0  & 0 \\\hline
         
 1& 1 & 1& \ldots  & 1 & 1 & 1 &  0& \ldots  & 0& 0 \\
  1& 1 & 1 & \ldots  & 1 & 1 & 1 &  1 & \ldots  & 0& 0 \\

    \vdots & \vdots &\vdots & \vdots & \vdots & \vdots & \vdots & \vdots & \vdots & \vdots & \vdots\\
    1&  1& 1 & \ldots  & 1 & 1 & 1 &  1 & \ldots  & 1& 0 \\
      1 & 1 & 1 & \ldots  & 1 & 1 & 1 &  1 & \ldots  & 1& 1 \\
         \end{array} \right]   \left[  \begin{array}{c}
x_{1}  \\
x_{2} \\
 \vdots \\
 x_{n-1} \\
 x_{n} \\\hline
 x_{n+1}\\
 x_{n+2}\\
 \vdots \\
 x_{2n-1} \\
x_{2n}\\
\end{array} \right]. $$

\noindent Since the lower triangular matrix involved has determinant $1$, it defines a element of ${\mathrm{SL}}_{2n}(\mathbb{Z})$ and therefore the 
 corresponding linear transfomation $f$ is an automorphism of the free $\mathbb{Z}$-module of rank $2n$. Thus, one can perform the bijective\footnote{It is crucial that the determinant of the transformation is equal to one.} substitution $\sigma \leftrightarrow x$ in the formula (\ref{defcnx}). Hence, one obtains
 
 \begin{equation}
 c_{n}=  \sum_{{\sigma \in \mathbb{Z}^{2n}}\atop{|\sigma_{i~\mathrm{mod}(n)}| \leqslant i}} \prod_{0 \leqslant i < k \leqslant n}\chi \left(1 \leqslant  (\sigma_{k}-\sigma_{i})^{2} + (\sigma_{k+n}-\sigma_{i+n})^{2} \leqslant( k-i)^{2} \right).
 \end{equation}

\noindent  Now, note the obvious but useful fact, saying that any positive integer greater than $1/2$ must be greater or equal than one. Hence, one can write 
 
 \begin{equation}
 c_{n}= \sum_{{\sigma \in \mathbb{Z}^{2n}}\atop{|\sigma_{i~\mathrm{mod}(n)}| \leqslant i}} \prod_{0 \leqslant i < k \leqslant n}\chi \left(1/2 \leqslant  (\sigma_{k}-\sigma_{i})^{2} + (\sigma_{k+n}-\sigma_{i+n})^{2} \leqslant( k-i)^{2} \right).
 \end{equation}
 The Lemma is proved.
 
 \begin{flushright}
 $\square$
 \end{flushright}

\noindent Given $\sigma \in \mathbb{Z}^{2n}$ and $0 \leqslant j < k \leqslant n$, let us define $$B_{jk}^{+}(\sigma) := B\left((\sigma_{j}, \sigma_{j+n}), k-j \right) \backslash B\left((\sigma_{j}, \sigma_{j+n}), \dfrac{1}{2} \right)$$ where $B(x,r)$ is the 2-dimensional Euclidean ball centered at $x\in \mathbb{R}^{2}$ and radius $r$. \\

\noindent The formula for $c_{n}$ may be written more concisely as

 \begin{equation}
 c_{n}= \sum_{\sigma \in \mathbb{Z}^{2n}}  \prod_{0 \leqslant l \leqslant n}     \chi_{[-l,l]^{2}}(\sigma_{l}, \sigma_{l+n}) \prod_{0 \leqslant j < k \leqslant n}\chi_{B^{+}_{jk}(\sigma)} \left( \sigma_{k}, \sigma_{k+n} \right).
 \end{equation}

\noindent A suitable way to write $c_{n}$ is as follows, 

\begin{equation}\label{uniformstyle}
 c_{n}= \sum_{\sigma \in \mathbb{Z}^{2n}}\prod_{0 \leqslant j < k \leqslant n} \left( \chi_{B^{+}_{jk}(\sigma)} \left( \sigma_{k}, \sigma_{k+n} \right) \prod_{l\neq j,k} \chi_{[-l,l]^{2}}(\sigma_{l}, \sigma_{l+n}) \right).
 \end{equation}

\noindent Indeed, if we denote by 
\begin{equation} \label{psi}
\psi_{jk}(\sigma) := \left( \chi_{B^{+}_{jk}(\sigma)} \left( \sigma_{k}, \sigma_{k+n} \right) \prod_{l\neq j,k} \chi_{[-l,l]^{2}}(\sigma_{l}, \sigma_{l+n}) \right)
\end{equation}

\noindent we get 
\begin{equation}
c_{n}= \sum_{\sigma \in \mathbb{Z}^{2n}}\prod_{0 \leqslant j < k \leqslant n} \psi_{jk}(\sigma).
\end{equation}
\subsection{An approximation of $c_{n}$} The functions $\psi_{jk}$ are not in the Schwartz class, but we can circonvent this issue by regularizing by a mollifier $(\eta_{\varepsilon})$ such that $\eta_{\varepsilon}$ weakly converges to $\delta$ as $\varepsilon \rightarrow 0^{+}$. Hence, the function $(\prod_{1 \leqslant j < k \leqslant n}\psi_{jk} )\ast \eta_{\varepsilon}$ is in the Schwartz class, thus, we can apply the Poisson summation formula

\begin{equation}\label{cncovol}
c_{n}=  \lim_{\epsilon \rightarrow 0^{+}} \sum_{\xi \in \mathbb{Z}^{2n}} \widehat{\left(\prod_{0 \leqslant j < k \leqslant n}\psi_{jk}  \ast \eta_{\varepsilon} \right)}(\xi).
\end{equation}

\noindent We know that the Fourier transform of a product is the convolution product of the Fourier transform hence, 
\begin{equation}\label{cncovol}
c_{n}=  \lim_{\epsilon \rightarrow 0^{+}} \sum_{\xi \in \mathbb{Z}^{2n}} \widehat{\left( \prod_{0 \leqslant j < k \leqslant n} \psi_{jk} \right)}(\xi) \widehat{\eta_{\varepsilon}}(\xi).
\end{equation}

\noindent The Fourier transform of a product of $L^{2}$ functions is the convolution product of such functions, thus 
\begin{equation}
c_{n}= \lim_{\epsilon \rightarrow 0^{+}}  \sum_{\xi \in \mathbb{Z}^{2n}} \left( \widehat{ \psi_{0,1}} \ast \widehat{ \psi_{0,2}} \ast \ldots  \ast \widehat{\psi_{n-1,n}}\right)(\xi) \widehat{\eta_{\varepsilon}}(\xi).
\end{equation}
Here the $m$-fold convolution of $m$ functions $f_{1}, \ldots, f_{m}$ on $\mathbb{R}^{2n}$ is given by 
$$\left( f_{1} \ast \ldots \ast f_{m} \right) (\xi)= \int_{\mathbb{R}^{2n}} \ldots \int_{\mathbb{R}^{2n}} f_{1}(x_{1}) \ldots f_{m-1}(x_{m-1}) f_{m}(\xi - \sum_{i=1}^{m-1} x_{i}) ~\bigotimes_{1\leqslant i \leqslant m-1} ~dx_{i}.$$

\noindent Thus, the number $c_{n}$ can be approximated by

\begin{equation}
c_{n}\sim  \sum_{\xi \in \mathbb{Z}^{2n}} \left( \widehat{ \psi_{0,1}} \ast \widehat{ \psi_{0,2}} \ast \ldots  \ast \widehat{\psi_{n-1,n}}\right)(\xi)
\end{equation}
since $ \widehat{\eta_{\varepsilon}}(\xi)$ tends to 1 weakly as $\varepsilon$ tends to 0.  

\subsection{The Fourier transform of $\psi_{jk}$.}

\noindent We are reduced to compute the Fourier transform of $\psi_{jk}$ for $0 \leqslant j < k \leqslant n$, it is given by the following proposition.
\begin{prop}\label{fourierprop} For every $\xi \in \mathbb{R}^{2n}$ and $0 \leqslant j<k \leqslant n$, the Fourier transform of $\psi_{jk}$ is given by 
\begin{equation}\label{fourierfunctions}
\begin{array}{c} \widehat{\psi_{jk}}(\xi) = \left(\displaystyle \prod_{{1\leqslant l \leqslant n}\atop{l \neq j,k}} 
l^{2}\sinc (2\pi l \xi_{l}) \sinc(2\pi l \xi_{l+n}) \right)  j^{2}\sinc(2\pi j (\xi_{j}+\xi_{k})) \sinc(2\pi j (\xi_{j+n}+\xi_{k+n})) \\
\left\lbrace \dfrac{(k-j) J_{1}(2\pi (k-j)\|(\xi_{k}, \xi_{k+n}) \|)  - \frac{1}{2} J_{1}(\pi \|(\xi_{k}, \xi_{k+n}) \|)  }{\|(\xi_{k}, \xi_{k+n}) \|} \right\rbrace,
\end{array}
\end{equation}
where $J_{1}$ is the Bessel function of the first kind of order 1. 
\end{prop}
\noindent \textbf{Proof.}  Let $\xi \in \mathbb{R}^{2n}$, the Fourier transform of $\psi_{jk}$ at $\xi$ is defined by 
$$ \widehat{\psi_{jk}} (\xi) =  \int_{\mathbb{R}^{2n}} \psi_{jk}(\sigma) e^{-2i\pi \sigma \cdot  \xi} ~ d\sigma.$$
Using the expression of $\psi_{jk}$ in (\ref{psi}), we get 
$$ \widehat{\psi_{jk}} (\xi) =  \int_{\mathbb{R}^{2n}} \chi_{B^{+}_{jk}(\sigma)} \left( \sigma_{k}, \sigma_{k+n} \right) \prod_{l\neq k} \chi_{[-l,l]^{2}}(\sigma_{l}, \sigma_{l+n}) e^{-2i\pi \sigma \cdot  \xi} ~ d\sigma.$$

$$ =  \int_{\mathbb{R}^{2n}} \chi_{B\left((\sigma_{j},\sigma_{j+n}),~ k-j \right)}  \left( \sigma_{k}, \sigma_{k+n} \right)  \prod_{l\neq k} \chi_{[-l,l]^{2}}(\sigma_{l}, \sigma_{l+n}) e^{-2i\pi \sigma \cdot  \xi} ~ d\sigma.$$
$$-  \int_{\mathbb{R}^{2n}} \chi_{B\left((\sigma_{j},\sigma_{j+n}),~ \frac{1}{2} \right)}  \left( \sigma_{k}, \sigma_{k+n} \right)  \prod_{l\neq k} \chi_{[-l,l]^{2}}(\sigma_{l}, \sigma_{l+n}) e^{-2i\pi \sigma \cdot  \xi} ~ d\sigma.$$

\noindent We need the following lemma, 
\begin{lem}\label{fourier} For any real $R>0$, one has 
$$ \int_{\mathbb{R}^{2n}} \chi_{B\left((\sigma_{j},\sigma_{j+n}),~ R \right)}  \left( \sigma_{k}, \sigma_{k+n} \right)  \prod_{l\neq k} \chi_{[-l,l]^{2}}(\sigma_{l}, \sigma_{l+n}) e^{-2i\pi \sigma \cdot  \xi} ~ d\sigma$$
$$ =  \left(\displaystyle \prod_{{1\leqslant l \leqslant n}\atop{l \neq j,k}} 
l^{2}\sinc (2\pi l \xi_{l}) \sinc(2\pi l \xi_{l+n}) \right) j^{2} \sinc(2\pi j (\xi_{j}+\xi_{l})) \sinc(2\pi j (\xi_{j+n}+\xi_{k+n})) \dfrac{ R J_{1}(2\pi R \|(\xi_{k}, \xi_{k+n}) \|) }{\|(\xi_{k}, \xi_{k+n}) \|}.$$
\end{lem}

\noindent \textit{Proof of the Lemma \ref{fourier}.}  By Fubini, we have the following, 
$$ \int_{\mathbb{R}^{2n}} \chi_{B\left((\sigma_{j},\sigma_{j+n}),~ R \right)}  \left( \sigma_{k}, \sigma_{k+n} \right)  \prod_{l\neq k} \chi_{[-l,l]^{2}}(\sigma_{l}, \sigma_{l+n}) e^{-2i\pi \sigma \cdot  \xi} ~ d\sigma$$

$$=  \int_{\mathbb{R}^{2n}} \chi_{B\left(0,~1 \right)}  \left( \frac{\sigma_{k}-\sigma_{j}}{R},  \frac{\sigma_{k+n}-\sigma_{j+n}}{R} \right)  \prod_{l\neq k} \chi_{[-l,l]^{2}}(\sigma_{l}, \sigma_{l+n}) e^{-2i\pi \sigma \cdot  \xi} ~ d\sigma$$

$$=  \int_{[-j,j]^{2}} e^{-2i\pi (\sigma_{j} \xi_{j}+ \sigma_{j+n}\xi_{j+n} )} d\sigma_{j} d\sigma_{j+n} \int_{\mathbb{R}^{2}}\chi_{B\left(0,~1 \right)}  \left( \frac{\sigma_{k}-\sigma_{j}}{R},  \frac{\sigma_{k+n}-\sigma_{j+n}}{R} \right) e^{-2i\pi (\sigma_{k} \xi_{k}+ \sigma_{k+n}\xi_{k+n} )}  d\sigma_{k} d\sigma_{k+n} $$

$$\left( \int \ldots \int_{\mathbb{R}^{2n-4}}\prod_{l\neq j,k} \chi_{[-l,l]^{2}}(\sigma_{l}, \sigma_{l+n}) e^{-2i\pi (\sigma_{l}  \xi_{l} + \sigma_{l+n}  \xi_{l+n})} ~ \bigotimes_{l\neq j,k} d\sigma_{l} \otimes d\sigma_{l+n} \right).$$

\noindent We perform the change of variables $$(\tau_{k},\tau_{k+n})= \varphi(\sigma_{k}, \sigma_{k+n}) := \left(\frac{\sigma_{k}-\sigma_{j}}{R},\frac{\sigma_{k}-\sigma_{j}}{R} \right). $$
\noindent  It is immediate to see that the jacobian of this substitution is equal to $R^{2}$. Thus, we get
$$ \int_{\mathbb{R}^{2n}} \chi_{B\left((\sigma_{j},\sigma_{j+n}),~ R \right)}  \left( \sigma_{k}, \sigma_{k+n} \right)  \prod_{l\neq k} \chi_{[-l,l]^{2}}(\sigma_{l}, \sigma_{l+n}) e^{-2i\pi \sigma \cdot  \xi} ~ d\sigma$$

$$=  \int_{[-j,j]^{2}} e^{-2i\pi (\sigma_{j} (\xi_{j}+ \xi_{j+n})+ \sigma_{j+n}(\xi_{j+n}+\xi_{k+n} ))} d\sigma_{j} d\sigma_{j+n} \int_{\mathbb{R}^{2}}\chi_{B\left(0,~1 \right)}  \left( \tau_{k}, \tau_{k+n} \right) e^{2i\pi (R\tau_{k} \xi_{k}+ R\tau_{k+n}\xi_{k+n} )}  R^{2}d\tau_{k} d\tau_{k+n} $$

$$\left( \int \ldots \int_{\mathbb{R}^{2n-4}}\prod_{l\neq j,k} \chi_{[-l,l]^{2}}(\sigma_{l}, \sigma_{l+n}) e^{-2i\pi (\sigma_{l}  \xi_{l} + \sigma_{l+n}  \xi_{l+n})} ~ \bigotimes_{l\neq j,k} d\sigma_{l} \otimes d\sigma_{l+n} \right).$$

\noindent Hence, 
$$\int_{\mathbb{R}^{2n}} \chi_{B\left((\sigma_{j},\sigma_{j+n}),~ R \right)}  \left( \sigma_{k}, \sigma_{k+n} \right)  \prod_{l\neq k} \chi_{[-l,l]^{2}}(\sigma_{l}, \sigma_{l+n}) e^{-2i\pi \sigma \cdot  \xi} ~ d\sigma$$ $$ = \widehat{\chi_{[-j,j]}}(\xi_{j}+\xi_{k})  \widehat{ \chi_{[-j,j]}}(\xi_{j+n}+ \xi_{k+n})   R^{2} \widehat{\chi_{B(0,1)}}(R\xi_{k}, R\xi_{k+n})   \prod_{l\neq j,k} \widehat{\chi_{[-l,l]}}(\xi_{l})  \widehat{ \chi_{[-l,l]}}(\xi_{l+n}) .$$

\noindent We need the two following well-known facts, \\

\begin{enumerate}
\item $ \widehat{ \chi_{[-l,l]}}(\xi) = 2l \sinc(2\pi l \xi).  $ for any $\xi \in \mathbb{R}$.
\item $ \widehat{ \chi_{B(0,1)}}(\xi_{1}, \xi_{2})= \dfrac{J_{1}(2\pi \|\xi \|)}{\| \xi \|}  $ where $J_{1}$ is the Bessel function of order 1 of the first kind and $\xi=(\xi_{1}, \xi_{2})\in \mathbb{R}^{2}$. (See Stein, \cite{stein})

\end{enumerate}
\noindent Using these statements we obtain the proof of Lemma \ref{fourier}.  
\begin{flushright}
$\square$
\end{flushright}

\noindent We are ready to compute the Fourier transform of $\psi_{jk}$. One has,

$$ \widehat{\psi_{jk}} (\xi)  =  \int_{\mathbb{R}^{2n}} \chi_{B\left((\sigma_{j},\sigma_{j+n}),~ k-j \right)}  \left( \sigma_{k}, \sigma_{k+n} \right)  \prod_{l\neq k} \chi_{[-l,l]^{2}}(\sigma_{l}, \sigma_{l+n}) e^{-2i\pi \sigma \cdot  \xi} ~ d\sigma$$
$$-  \int_{\mathbb{R}^{2n}} \chi_{B\left((\sigma_{j},\sigma_{j+n}),~ \frac{1}{2} \right)}  \left( \sigma_{k}, \sigma_{k+n} \right)  \prod_{l\neq k} \chi_{[-l,l]^{2}}(\sigma_{l}, \sigma_{l+n}) e^{-2i\pi \sigma \cdot  \xi} ~ d\sigma.$$

\noindent Thus, using Lemma \ref{fourier} for both integral taking $R=k-j$  for the first integral and $R=1/2$ for the second we get the Fourier transform as in (\ref{fourierfunctions}). The proposition \ref{fourierprop} is proved. 

\begin{flushright}
$\square$
\end{flushright}

\subsection{An integral asymptotic formula for $c_n$ and end of the proof of Theorem \ref{main}}

\noindent We derive from (\ref{cncovol}), the following expression for $c_{n}$,
\begin{equation}
c_{n} = \lim_{\epsilon \rightarrow 0^{+}}\sum_{v \in \mathbb{Z}^{2n}} \left(\widehat{\psi_{0,1}} \ast \ldots \ast \widehat{\psi_{n-1,n}} \right) \left( v \right) \widehat{\eta_{\varepsilon}}(v).
\end{equation}

\noindent Hence, using Proposition \ref{fourierprop} we obtain the following estimate in the weak sense  as $n \rightarrow \infty$,  
\begin{equation}~\label{form}
c_{n} \sim \sum_{v \in \mathbb{Z}^{2n}} \int_{\mathbb{R}^{2n}} \ldots \int_{\mathbb{R}^{2n}} \widehat{\psi_{0,1}}(\xi_{0,1}) \ldots \widehat{\psi_{n-2,n}}(\xi_{n-2,n})$$ $$ \widehat{\psi_{n-1,n}}\left( v -(\sum_{1\leqslant j<k\leqslant n-1 }\xi_{j,k} + \sum_{1 \leqslant j \leqslant n-2} \xi_{j,n})\right) \left(\bigotimes_{0 \leqslant j<k\leqslant n-1} d\xi_{j,k}\right) \bigotimes \left(\bigotimes_{0\leqslant j \leqslant n-2} d\xi_{j,n} \right)
 \end{equation}
 where \begin{equation}\label{fourierfunctions}
\begin{array}{c} \widehat{\psi_{jk}}(\xi) = \left(\displaystyle \prod_{{1\leqslant l \leqslant n}\atop{l \neq j,k}} 
l^{2}\sinc (2\pi l \xi_{l}) \sinc(2\pi l \xi_{l+n}) \right)  j^{2}\sinc(2\pi j (\xi_{j}+\xi_{k})) \sinc(2\pi j (\xi_{j+n}+\xi_{k+n})) \\
\left\lbrace \dfrac{(k-j) J_{1}(2\pi (k-j)\|(\xi_{k}, \xi_{k+n}) \|)  - \frac{1}{2} J_{1}(\pi \|(\xi_{k}, \xi_{k+n}) \|)  }{\|(\xi_{k}, \xi_{k+n}) \|} \right\rbrace.
\end{array}
\end{equation}

\noindent Each variable $\xi_{j,k}$ is a vector of $\mathbb{R}^{2n}$ whose components are given by 
$$ \xi_{j,k} = \left( \xi_{j,k}^{1}, \ldots,    \xi_{j,k}^{2n}\right).$$

\noindent The multivariate integral \ref{form} above involves $\frac{(n-2)(n-1)}{2} + n-2 = \frac{(n-2)(n+1)}{2}$ vectors in $\mathbb{R}^{2n}$. Thus, one can see this integral as an integral in $n(n-2)(n+1)$ real variables. 

\noindent Finally, we obtain the required estimate as $n$ gets large,
$$ c_{n} \sim \sum_{v \in \mathbb{Z}^{2n}} \int_{\mathbb{R}^{n(n-2)(n+1)}} \prod_{0 \leqslant j<k\leqslant n-1} \widehat{\psi_{j,k}}(\xi_{j,k}) $$ $$ \prod_{0 \leqslant j \leqslant n-2 }  \widehat{\psi_{j,n}}(\xi_{j,n}) ~ \widehat{\psi_{n-1,n}}\left( v -(\sum_{1\leqslant j<k\leqslant n-1 }\xi_{j,k} + \sum_{1 \leqslant j \leqslant n-2} \xi_{j,n}) \right) \left(\bigotimes_{0 \leqslant j<k\leqslant n-1} d\xi_{j,k}\right) \bigotimes \left(\bigotimes_{0\leqslant j \leqslant n-2} d\xi_{j,n} \right). $$\\
 where
$$\begin{array}{c} \widehat{\psi_{jk}}(\xi) = \left(\displaystyle \prod_{{1\leqslant l \leqslant n}\atop{l \neq j,k}} 
l^{2}\sinc (2\pi l \xi_{l}) \sinc(2\pi l \xi_{l+n}) \right)  j^{2}\sinc(2\pi j (\xi_{j}+\xi_{k})) \sinc(2\pi j (\xi_{j+n}+\xi_{k+n})) \\
\left\lbrace \dfrac{(k-j) J_{1}(2\pi (k-j)\|(\xi_{k}, \xi_{k+n}) \|)  - \frac{1}{2} J_{1}(\pi \|(\xi_{k}, \xi_{k+n}) \|)  }{\|(\xi_{k}, \xi_{k+n}) \|} \right\rbrace.\end{array}
$$

\noindent This achieves the proof of Theorem \ref{main}. 

\section{How to obtain an exact Asymptotic estimate of $c_{n}$?}

\subsection{The main term of the asymptotic expansion for $c_{n}$} The main theorem gives us an asymptotic involving an integral of a product of the transcendental functions $\sinc$ and $J_{1}$. For the Poisson sum, the volume term corresponding to $v=0$ is the main term since the Fourier coefficients decreases as $\| v\|$ tends to infinity.\\
\noindent A first remark is that this the functions$\sinc$ and Bessel $J_{1}$ have their mass concentrated near the zero. Therefore, most of the contribution of the integral in Theorem \ref{main} will be concentrated near the zero. We will consider the truncated integral for $|\xi|_{\infty}$ less than the least first zero of both $J_{1}(x)/x$ and $\sinc(x)$. Indeed, we can assume that we are integrating a product of positive functions which are concentrated near the zero, both have quite similar behaviours, in other words, they will not neutralize each other. 
The similarity of both functions is illustrated by their infinite product expansion which have the same nature, namely, 

\begin{equation}\label{sincprodformula} \sinc(2\pi x ) = \prod_{m\geqslant 1} \left(1-\dfrac{4x^2}{m^{2}} \right)
\end{equation}

\noindent and 

\begin{equation}\label{besselprodformula} \dfrac{J_{1}(x)}{x} = \dfrac{1}{2}  \prod_{m\geqslant 1} \left(1-\dfrac{x^2}{j_m^{2}} \right)
\end{equation}
where $j_{1}< j_{2}< \ldots$ is the set of zeros of  the Bessel function $J_{1}(x)$.

\subsection{Some notations} Let us set  $\displaystyle U_{n}(\xi)= \prod_{1 \leqslant j<k\leqslant n-1} \widehat{\psi_{j,k}}(\xi_{j,k})$,  $\displaystyle V_{n}(\xi)= \prod_{1 \leqslant j \leqslant n-2} \widehat{\psi_{j,n}}(\xi_{j,n})$ and 
$$ W_{n}(\xi,v) = \widehat{\psi_{n-1,n}}\left( v -(\sum_{1\leqslant j<k\leqslant n-1 }\xi_{j,k} + \sum_{1 \leqslant j \leqslant n-2} \xi_{j,n}) \right)$$
for $\xi= (\xi_{1,2}^{1}, \ldots, \xi_{1,2}^{2n}, \ldots, \xi_{n-2,n}^{1}, \ldots, \xi_{n-2,n}^{2n}) \in \mathbb{R}^{n(n-2)(n+1)}$. Therefore, as $n \rightarrow \infty$ one has
\begin{equation}
c_{n} \sim \sum_{v\in \mathbb{Z}^{2n-2}} I_{n}(v)
\end{equation}
where 
$$I_{n}(v) :=  \int_{\mathbb{R}^{n(n-2)(n+1)}} U_{n}(\xi) V_{n}(\xi) W_{n}(\xi,v) ~d\xi.$$
 
\noindent We have seen in Proposition \ref{fourierprop}, that the Fourier transform of $\psi_{jk}$ can be expressed as follows
\begin{equation}
\begin{array}{c} \widehat{\psi_{jk}}(\xi_{jk}) = \dfrac{(n!)^{2}}{k^{2}}\left(\displaystyle \prod_{{1\leqslant l \leqslant n}\atop{l \neq j,k}} 
\sinc (2\pi l \xi_{jk}^{l}) \sinc(2\pi l \xi_{jk}^{l+n}) \right) \sinc(2\pi j (\xi_{jk}^{j}+\xi_{jk}^{k})) \sinc(2\pi j (\xi_{jk}^{j+n}+\xi_{jk}^{k+n})) \\
\left\lbrace \dfrac{(k-j) J_{1}(2\pi (k-j)\|(\xi_{jk}^{k}, \xi_{jk}^{k+n}) \|)  - \frac{1}{2} J_{1}(\pi \|(\xi^{k}_{jk}, \xi^{k+n}_{jk}) \|)  }{\|(\xi^{k}_{jk}, \xi^{k+n}_{jk}) \|}\right\rbrace.
\end{array}
\end{equation}

\noindent Thus, 
\begin{equation}
\begin{array}{c}
U_{n}(\xi)= \prod_{1 \leqslant j<k \leqslant n} \dfrac{(n!)^{2}}{k^{2}}\left(\displaystyle \prod_{{1\leqslant l \leqslant n}\atop{l \neq j,k}} 
\sinc (2\pi l \xi_{jk}^{l}) \sinc(2\pi l \xi_{jk}^{l+n}) \right) \sinc(2\pi j (\xi_{jk}^{j}+\xi_{jk}^{k})) \sinc(2\pi j (\xi_{jk}^{j+n}+\xi_{jk}^{k+n})) \\
\left\lbrace \dfrac{(k-j) J_{1}(2\pi (k-j)\|(\xi_{jk}^{k}, \xi_{jk}^{k+n}) \|)  - \frac{1}{2} J_{1}(\pi \|(\xi^{k}_{jk}, \xi^{k+n}_{jk}) \|)  }{\|(\xi^{k}_{jk}, \xi^{k+n}_{jk}) \|}\right\rbrace.
\end{array}
\end{equation}

\noindent and 

\begin{equation}
\begin{array}{c}
V_{n}(\xi)= \prod_{1 \leqslant j \leqslant n-2} \dfrac{(n!)^{2}}{n^{2}}\left(\displaystyle \prod_{{1\leqslant l \leqslant n}\atop{l \neq j,n}} 
\sinc (2\pi l \xi_{jn}^{l}) \sinc(2\pi l \xi_{jn}^{l+n}) \right) \sinc(2\pi j (\xi_{jn}^{j}+\xi_{jn}^{k})) \sinc(2\pi j (\xi_{jn}^{j+n}+\xi_{jn}^{k+n})) \\
\left\lbrace \dfrac{(k-j) J_{1}(2\pi (k-j)\|(\xi_{jk}^{k}, \xi_{jk}^{k+n}) \|)  - \frac{1}{2} J_{1}(\pi \|(\xi^{k}_{jk}, \xi^{k+n}_{jk}) \|)  }{\|(\xi^{k}_{jk}, \xi^{k+n}_{jk}) \|} \right\rbrace.
\end{array}
\end{equation}

\subsection{The main contribution}

Even though a closed formula in terms of classical functions seems out of reach, we can try to get an asymptotic formula by identifying the main contribution of the oscillatory integral $I_{n}(v)$ as $n \rightarrow \infty$. Since $I_{n}(v)$ is rapidly decreasing to zero as a the Fourier transform of $L^{1}$ function, the major contribution of $c_{n}$ is given by the value $I_{n}(0)$. In other words, as $n$ gets large 
$$ c_{n} \sim I_{n}(0) \left( 1+ o(1) \right).$$
 
 \noindent The integral $I_{n}(0)$ consists in a product of two kind of functions namely$\sinc$ and Bessel $J_{1}$ functions. Those two functions satisfies the following asymptotics at when $x$ tends to zero 
 \begin{equation}
 \sinc(x)= O(1) \ \ \ \ \  \ \mathrm{and} \ \ \ \ \ \ \frac{J_{1}(|x|)}{x}= O(1).
 \end{equation}
 
\noindent  Moreover, $\sinc$ and $J_{1}(x)/x$ achieve their maximum at $x=0$, hence the Fourier factors in $I_{n}(0)$ are maximized when $\xi =0 $ (see figure 1).  Hence, the mass of the integral $I_n = I_{n}(0)$ is essentially concentrated around zero. Thus, if we split the integral  $I_{n}$ as follows
 \begin{equation}
 I_{n} = \int_{\|\xi\|_{\infty} \leqslant \delta}  U_{n}(\xi) V_{n}(\xi) W_{n}(\xi,0) ~d\xi+   \int_{\|\xi\|_{\infty} \geqslant \delta}U_{n}(\xi) V_{n}(\xi) W_{n}(\xi,0) ~d\xi
 \end{equation}
 for some $\delta > 0$ small enough. If the first integral $I_{n}^{\delta}$ dominates the second integral $R_{n}^{\delta}$,  we will have 
 \begin{equation}
 c_{n} \sim I_{n}^{\delta} \left(1 + o(1)\right).
 \end{equation}
 
 Our hope is that as for a well-chosen $\delta$, $I_{n}^{\delta}$ would follow the expected asymptotic. 
 
\section{The mean squared displacement of the endpoint of a typical SAW}

 \noindent Let $\omega_{n}$ be randomly chosen self-avoiding of lenght $n$ in the square lattice, we define the \textit{ mean squared displacement}  of $\omega_{n}$ to its endpoint by the following quantity
$$\langle\| w_{n} \|^{2}  \rangle=  \mathbf{E}(\| w_{n} \|^{2}) = \dfrac{1}{| \Omega_{n}|} \sum_{\omega \in \Omega_{n}} \| \omega_{n} \|^{2} = \dfrac{1}{ c_{n}} \sum_{\omega \in \Omega_{n}} \| \omega_{n} \|^{2}.$$

\noindent For a simple random walk of lentgh $n$ in the square lattice, the mean squared displacement is of order $\sqrt{n}$.  Therefore, it is natural to expect that for a self-avoiding walk, its endpoint tends to move away from the origin in order to avoid self-intersections. The prediction of Paul Flory  \cite{flory} is that $$\langle\| w_{n} \|^{2}  \rangle \sim n^{3/2+o(1)}.$$ 
\noindent The best current bound is due to Duminil-Copin and Hammond  \cite{dmsub} who showed that 
$$\langle\| w_{n} \|^{2} \rangle=o(n^{2}).$$
 The hope is to improve the exponent $\gamma=2$  in order to reach  Flory's prediction.  We proceed as for $c_{n}$, in that, we interpret the sum $\sum_{\omega \in \Omega_{n}} \| \omega_{n} \|^{2}$ as a counting function. Using the same notations as \S 2, we can write 
 $$\sum_{\omega \in \Omega_{n}} \| \omega_{n} \|^{2} = \sum_{ w \in \mathbb{Z}^{2n}}  \|w \|^{2} \prod_{0 \leqslant j < k \leqslant n} \chi \left(1/2 \leqslant \| w_{k}-w_{j} \|^{2} \leqslant (k-j)^{2} \right).$$
 
 $$=\sum_{{x\in \mathbb{Z}^{2n}}\atop{\|x\|_{\infty} \leqslant 1}} \left[ (x_{1}+ \ldots+ x_{n})^{2}+  (x_{n+1}+ \ldots+ x_{2n})^{2} \right] \prod_{0 \leqslant j < k \leqslant n} \chi \left(1/2 \leqslant (x_{j+1}+ \ldots+ x_{k})^{2}+  (x_{n+j+1}+ \ldots+ x_{n+k})^{2} \leqslant (k-j)^{2} \right).$$
 
\noindent We make the bijective substitution $x \leftrightarrow \sigma$ of Lemma \ref{cnsigma}, therefore we obtain 
$$\sum_{\omega \in \Omega_{n}} \| \omega_{n} \|^{2}=\sum_{{\sigma \in \mathbb{Z}^{2n}}\atop{\|\sigma_{l}\|_{\infty} \leqslant l}} \left[ \sigma_{n}^{2}+  \sigma_{2n}^{2} \right]  \prod_{0 \leqslant j < k \leqslant n} \chi \left(1/2 \leqslant (\sigma_{k}-\sigma_{j})^{2}+  (\sigma_{n+k}-\sigma_{n+j})^{2} \leqslant (k-j)^{2} \right).$$

\noindent Again using the notations of section \S 2, we have
$$\sum_{\omega \in \Omega_{n}} \| \omega_{n} \|^{2}=\sum_{{\sigma \in \mathbb{Z}^{2n}}} \left[ \sigma_{n}^{2}+  \sigma_{2n}^{2} \right]  \prod_{0 \leqslant j < k \leqslant n} \left( \chi_{B^{+}_{jk}(\sigma)} \left( \sigma_{k}, \sigma_{k+n} \right) \prod_{l\neq j,k} \chi_{[-l,l]^{2}}(\sigma_{l}, \sigma_{l+n}) \right).$$

\noindent Hence, 

$$\sum_{\omega \in \Omega_{n}} \| \omega_{n} \|^{2}= \sum_{{\sigma \in \mathbb{Z}^{2n}}} \left[ \sigma_{n}^{2}+  \sigma_{2n}^{2} \right]  \chi_{[-n,n]^{2}}(\sigma_{n}, \sigma_{2n})\prod_{0 \leqslant j < k \leqslant n} \psi_{jk}(\sigma).$$

\noindent Then, if we set $ \psi_{n}(\sigma) =  \left[ \sigma_{n}^{2}+  \sigma_{2n}^{2} \right] \prod_{1 \leqslant l \leqslant n} \chi_{[-l,l]^{2}} (\sigma_{l}, \sigma_{l+n})$

$$\sum_{\omega \in \Omega_{n}} \| \omega_{n} \|^{2} = \sum_{{\sigma \in \mathbb{Z}^{2n}}} \psi_{n}(\sigma) \prod_{0 \leqslant j < k \leqslant n} \psi_{jk}(\sigma).$$

\noindent After regularizing by a mollifier and applying Poisson summation formula, we get

$$\sum_{\omega \in \Omega_{n}} \| \omega_{n} \|^{2} \sim \sum_{{\xi \in \mathbb{Z}^{2n}}} \left( \widehat{\psi_{n}} \ast  \widehat{\psi_{0,1}} \ast \ldots \ast \widehat{\psi_{n-1,n}} \right)(\xi) $$
\noindent where $\widehat{\psi_{j,k}}$ have been computed in Proposition  \ref{fourierprop} and for $\xi=(\xi_{1}, \ldots , \xi_{2n})$
 $$  \widehat{\psi_{n}}(\xi) = \prod_{1 \leqslant l \leqslant n-1} \sinc(\pi l (\xi_{l}+ \xi_{l+n}))  \left( \dfrac{\partial^{2}}{\partial \xi_{n}^{2}}[ \sinc(\pi n \xi_{n})] + \dfrac{\partial^{2}}{\partial \xi_{2n}^{2}}[ \sinc(\pi n \xi_{2n}) ] \right).$$
Thus, 
 $$\langle\| w_{n} \|^{2}  \rangle \sim \dfrac{\sum_{{\xi \in \mathbb{Z}^{2n}}} \left( \widehat{\psi_{n}} \ast  \widehat{\psi_{0,1}} \ast \ldots \ast \widehat{\psi_{n-1,n}} \right)(\xi)}{c_{n}} \sim  \dfrac{\sum_{{\xi \in \mathbb{Z}^{2n}}} \left( \widehat{\psi_{n}} \ast  \widehat{\psi_{0,1}} \ast \ldots \ast \widehat{\psi_{n-1,n}} \right)(\xi)}{\sum_{{\xi \in \mathbb{Z}^{2n}}} \left(   \widehat{\psi_{0,1}} \ast \ldots \ast \widehat{\psi_{n-1,n}} \right)(\xi)}. $$
 
\section{Final remarks}

\noindent In conclusion, we have adressed a new way to looking at the problem of enumeration of self-avoiding walks. The method consists in translating the problem into a lattice points counting in multidimensional domains. It gives concrete asymptotic formulas for $c_{n}$ in terms of classical transcendental functions and a unified set-up which allows to study  SAWs of any dimension. To sum-up the lattice point counting method gives a natural and unified approach to the study of the asymptotic behaviour for
\begin{enumerate}
\item the number of SAWs in $\mathbb{Z}^{2}$, $c_{n}$  in dimension 2.
\item  The number of SAWs in $\mathbb{Z}^{d}$, $c_{n}(d)$ for $d \geq 3$.
\item The mean squared displacement $\langle \|w\|^{2} \rangle$ in dimension $2$.
\end{enumerate}

The multidimensional integral involved in the main theorem is, essentially, of \textit{Gaussian type}. In order words, we are reduced to compute determinants of hypermatrices. We hope that in the near future, we will be able to estimate such integrals given any dimension.  This would provide accurate bounds for the number of self-avoiding walks of a given length and for the expected mean displacement for such walks.


\begin{thebibliography}{ZZZZ}
 \bibitem{bauer}  BAUERSCHMIDT, R., DUMINIL-COPIN, H., GOODMAN, J. and SLADE, G. (2012). Lectures on self-
avoiding walks. In Probability and Statistical Physics in Two and More Dimensions. Clay Math. Proc.
15 395–467. Amer. Math. Soc., Providence.
 \bibitem{beaton} BEATON, N. R., BOUSQUET-MÉLOU, M., DE GIER, J., DUMINIL-COPIN, H. and
GUTTMANN, A. J. (2014). The critical fugacity for surface adsorption of self-avoiding
walks on the honeycomb lattice is $1+\sqrt{2}$ Comm. Math. Phys. 326 727–754.

 \bibitem{brydges}  BRYDGES, D. C., IMBRIE, J. Z. and SLADE, G. (2009). Functional integral representations
for self-avoiding walk. Probab. Surv. 6 34–61

\bibitem{flory}  FLORY P.J.  The configuration of a real polymer chain, J. Chem. Phys. 17 (1949), 303–310.



 \bibitem{dmsub}  DUMINIL-COPIN, H. and HAMMOND, A. (2013). Self-avoiding walk is sub-ballistic. Comm.
Math. Phys. 324 401–423

 \bibitem{dmhex}    DUMINIL-COPIN H. and  SMIRNOV S. The connective constant of the honeycomb lattice equals $\sqrt{2+\sqrt{2}}$, Annals of Math. 175(3), 1653-1665 (2012).
 \bibitem{conway}   CONWAY A.R.,  ENTING I.G.  and  GUTTMAN A.J. Algebraic techniques for enumerating self-
avoiding walks on the square lattice, J. Phys. A: Math. Gen. 26 (1993), 1519–1534. 

 \bibitem{glazman} GLAZMAN, A. Connective constant for a weighted self-avoiding walk on $Z^{2}$. Electron. Commun. Probab. 20: 1-13 (2015).
 
 \bibitem{hara}  HARA  T.  and SLADE G. , Self-avoiding walk in five or more dimensions. I. The critical behaviour, Commun. Math. Phys. 147 (1992), 101–136.



\bibitem{hw}   HAMMERSLEY  J.M. and WELSH  D.J.A., Further results on the rate of convergence to the con
nective constant of the hypercubical lattice, Quart. J. Math. Oxford (2), 13 (1962), 108–110.

 \bibitem{kesten} KESTEN, H. (1963). On the number of self-avoiding walks. J. Math. Phys. 4 960–969.

 \bibitem{madras1} MADRAS, N. (2014). A lower bound for the end-to-end distance of the self-avoiding walk.
Canad. Math. Bull. 57 113–118.

\bibitem{madras-slade}  Madras N. and Slade G. The self-avoiding walk, Birkh\"{a}user, Boston, (1993).
\bibitem{nienhuis} NIENHUIS B.  Exact critical exponents of the $O(n)-$models in two dimensions, Phys. Rev.
Lett. 49 (1982), 1062–1065.


\bibitem{sle}  LAWLER G.F.,  SCHRAMM O.  and  WERNER W. On the scaling limit of planar self-avoiding walk.
Proc. Symposia Pure Math., 72:339–364, (2004).

 \bibitem{orr} ORR, W. J. C. (1947). Statistical treatment of polymer solutions at infinite dilution. Trans.
Faraday Soc. 43 12–27.


\bibitem{stein}STEIN  E. M. , Harmonic Analysis: Real-Variable Methods, Orthogonality, and Oscillatory Integrals, Princeton University Press; First Edition (July 12, 1993)
\end{thebibliography}
\end{document}